\documentclass{svmult}
\usepackage{makeidx}         
\usepackage{graphicx}        
\usepackage{multicol}        
\usepackage[bottom]{footmisc}
\usepackage{amsmath,amssymb,euscript,psfrag,latexsym,graphicx}

\makeindex             
                      
\begin{document}
\newcommand{\mR}{{\mathbb R}}
\newcommand{\mZ}{{\mathbb Z}}
\newcommand{\mC}{{\mathbb C}}
\newcommand{\fR}{{\mathfrak R}}
\newcommand{\fG}{{\mathfrak G}}
\newcommand{\fT}{{\mathfrak T}}
\newcommand{\fC}{{\mathfrak C}}
\newcommand{\fF}{{\mathfrak F}}
\newcommand{\fM}{{\mathfrak M}}

\newcommand{\bR}{{\bf R}}
\newcommand{\hR}{\hat{\bf R}}
\newcommand{\tR}{\tilde{\bf R}}
\newcommand{\bQ}{{\bf Q}}
\newcommand{\hQ}{\hat{\bf Q}}
\newcommand{\tQ}{\tilde{\bf Q}}
\newcommand{\bX}{{\bf X}}
\newcommand{\cM}{{\mathcal M}_n}
\newcommand{\cT}{{{\mathcal T}_n}}
\newcommand{\cTi}{{_{\mathcal T}}}
\newcommand{\cC}{{\mathcal C}}
\newcommand{\trace}{{\rm trace\,}}

\newcommand{\by}{{\bf y}}\newcommand{\bu}{{\bf u}}

\newcommand{\hy}{{\hat{\bf y}}}
\newcommand{\bw}{{{\bf w}}}
\newcommand{\bz}{{{\bf z}}}
\newcommand{\ty}{{\tilde{\bf y}}}
\newcommand{\intpi}{\frac{1}{2\pi}\int_{-\pi}^\pi}
\newcommand{\rR}{{\rm R}}

\title*{Distances between time-series and their autocorrelation statistics}
\titlerunning{Distances between time-series \& statistics}
\author{Tryphon T. Georgiou}
\institute{Department of Electrical and Computer Engineering, University of Minnesota\\ Minneapolis, MN 55455; email: \texttt{tryphon@umn.edu}\\
$\,$\\
{\em This paper is dedicated to Giorgio Picci on the occasion of his sixty-fifth birthday}
}
\maketitle

\begin{abstract}
We begin with an interpretation of the $L_1$-distance between two power spectral densities and then, following an analogous rationale, we develop a natural metric for
quantifying distance between respective covariance matrices.
\end{abstract}

\section{Introduction}

Consider two discrete-time, stationary, zero-mean, (real-valued for notational convenience) random processes $\by_k$ and $\hy_k$ ($k\in\mZ$) having power spectral densities $f_\by(\theta)$ and $f_\hy(\theta)$ ($\theta\in[-\pi,\pi]$), and autocorrelation functions $R_\ell$ and $\hat{R}_\ell$ ($\ell\in\mZ$), respectively, i.e.,
\[ R_\ell = E\{\by_k\by_{k+\ell}\} = \intpi f(\theta) e^{-j\ell \theta} d\theta,
\]
and similarly for the ``hatted'' quantities. When the power spectrum contains a singular part, then $f(\theta)d\theta$ needs to be replaced by a non-negative finite spectral measure $d\mu(\theta)$.   

We are interested in quantifying the distance between respective spectra and statistics for two such random process $\by_k$ and $\hy_k$. When two vectors
\begin{eqnarray*}
\rR_n&:=&\left[\begin{matrix}R_0&R_1&\ldots &R_{n-1}\end{matrix}\right], \mbox{ and}\\
\hat\rR_n&:=&\left[\begin{matrix}R_0&R_1&\ldots &R_{n-1}\end{matrix}\right]
\end{eqnarray*}
of autocorrelation samples
are available and need to be compared, one may use any metric in $\mR^{n}$ for that purpose, as for instance $\|\rR_n-\hat\rR_n\|_2=\sqrt{\sum_k (R_k-\hat R_k)^2}$. However, we are not aware of any significance that can be attached to such a distance other than the fact that it is a metric in $\mR^{n}$. Our goal in this paper, is to seek a metric which can be physically motivated.

Similarly, if we are to compare $f_\by(\theta)$ and $f_\hy(\theta)$, it appears difficult to motivate the use of an $L_2$-distance $\|f_\by(\theta)-f_\hy(\theta)\|_2$. For one thing, the $L_2$-distance cannot be generalized to deal with spectral measures when singular parts are present. There are certainly other alternatives. In the speech processing literature in particular there is a plethora of distances that, however, are not metrics \cite{Gray} but have been motivated by specific needs. Function theoretic alternatives that one can use (e.g., $L_p$-norms, etc.) including Wasserstein-like transportation measures typically lack a physical interpretation. In a recent study \cite{Geo2007a,Geo2007b} a pseudometric was constructed as a geodesic between spectral densities/measures with respect to a rather natural Riemannian metric ---this metric quantifies the degradation of predictive-error variance when the predictor is designed based on the wrong choice between two alternatives and the geometry is, in essence, Euclidean but only after we transform spectral densities using the logarithmic map.

In the current paper we focus on the $L_1$ distance
\[
\|f_\by(\theta)-f_\hy(\theta)\|_1:=\intpi |f_\by(\theta)-f_\hy(\theta)|\,d\theta
\]
which has also a rather natural interpretation. After a brief discussion of the relevance of the $L_1$-distance, following a similar rationale, we will develop an analogous metric between finite partial covariance data of the corresponding random processes.

\section{Interpretation of the $L_1$ distance}

Given $\by_k$ and $\hy_k$ we postulate that ther exist
two random processes $\psi_k$ and $\hat{\psi}_k$ so that
\begin{equation}\label{equalityholds}
\by_k+\psi_k=\hy_k+\hat{\psi}_k.
\end{equation}
Alternatively, we postulate that there exists a random process $\bz_k$ and that the two original random processes relate to $\bz_k$ via
\begin{eqnarray*}
\by_k&=&\bz_k-\psi_k\\
\hy_k&=&\bz_k-\hat\psi_k.
\end{eqnarray*}
It is natural to seek such perturbations
of minimal total combined variance
\begin{equation}\label{combined} E\{\psi_k^2\}+E\{\hat{\psi}_k^2\}
\end{equation}
that is sufficient to ``reconcile'' the two processes.
The combined variance $E\{\psi_k^2\}+E\{\hat{\psi}_k^2\}$ represents the minimal amount of ``energy'' of perturbations in the two time-series that is needed to render the two indistinguishable. Intuitively, the minimal combined variance which is consistent with the available data quantifies the distance between the two. 

Given $f_\by,f_\hy$, the optimal choice consists of random processes $\psi_k$ and $\hat\psi_k$ such that $\psi_k$ and $\hat\psi_k$ are independent, $\hy_k$ and $\hat{\psi}_k$ are also independent, and
\begin{subequations}
\begin{eqnarray}\label{psi}
f_\psi(\theta) =
\left\{\begin{array}{cl} f_\hy(\theta)-f_\by(\theta)&\mbox{ if }f_\hy(\theta)-f_\by(\theta)\geq 0,\\
                                 0 &\mbox{ otherwise,}\end{array}\right.\\\label{psihat}
f_{\hat{\psi}}(\theta) =\left\{\begin{array}{cl}f_\by(\theta) -f_\hy(\theta)&\mbox{ if }f_\hy(\theta)-f_\by(\theta)\leq 0,\\ 0 &\mbox{ otherwise.}\end{array}\right.
\end{eqnarray}
\end{subequations}
Then, the power spectrum of the ``sum''
\begin{equation}\nonumber
\bz_k:=\by_k+\psi_k=\hy_k+\hat{\psi}_k
\end{equation}
is simply
\[f_z(\theta):=\max\{f_\by(\theta),\,f_\hy(\theta)\},\;\theta\in[-\pi,\pi],
\]
and
\begin{eqnarray}\label{metric}
d(f_\by,f_\hy)&:=& E\{\psi_k^2\}+E\{\hat{\psi}_k^2\}\\\nonumber
&=&\intpi (f_\psi(\theta)+f_{\hat{\psi}}(\theta))d\theta\\\nonumber
&=&\intpi | f_\by(\theta)-f_\hy(\theta)|\,d\theta\\\label{L1}
&=& \|f_\by-f_\hy\|_1.
\end{eqnarray}
Obviously, this construction extends in the obvious way to the case of not-necessarily absolutely continuous power spectra as well, and the metric includes the measure of any discrepancies between the singular parts of the two spectral measures.\footnote{It will be interesting to explore the practical significance of other possibilities for quantifying distance such as
\[
\frac{\int^\pi_{-\pi} (f_\psi(\theta)+f_{\hat{\psi}}(\theta))d\theta}{\int^\pi_{-\pi}f_\bz(\theta) d\theta}
\mbox{ or }
\int^\pi_{-\pi} \frac{  f_\psi(\theta)+f_{\hat{\psi}}(\theta)} {f_\bz(\theta)} \frac{d\theta}{2\pi}.
\]}
Clearly, $d(f_\by,f_\hy)$ is a metric as seen from (\ref{L1}). Building on a similar rationale, in the next section, we develop a metric for covariance matrices.

\section{A distance between covariance matrices}

It is often the case that only a finite segment of the autocorrelation function of time-series $\by_k$ and $\hy_k$ is available (and even then, possibly uncertain). Thus, it is of interest to consider distances between the partial autocorrelation statistics $\rR$ and
$\hat\rR$.
To this end, we follow the dictum of the previous section and define as a distance measure
the minimal combined variance of random processes $\psi_k$ and $\hat{\psi}_k$ for which (\ref{equalityholds}) holds. Naturally, since only partial covariance samples are available,  the random processes $\psi_k,\hat{\psi}_k$ and (\ref{equalityholds}) need to be consistent with these data.

First denote by
\[\bR_n:=\left[\begin{array}{cccc}R_0&R_1&\ldots&R_{n-1}\\
R_{-1}&R_0& \ldots& R_{n-2}\\
\vdots & \vdots& \ddots & \vdots\\
R_{-(n-1)}& R_{-(n-2)}&\ldots &R_0\end{array}\right]
\]
the $n\times n$ covariance matrix corresponding to $\by_k$ and the covariance samples in $\rR_n$ and, similarly, $\hat{\bR}_n$ for the Toeplitz matrix based on $\hat{\rR}_n$.
If $\bQ_n,\hQ_n$ denote the corresponding finite Toeplitz covariances of the random processes $\psi_k$ and $\hat{\psi}_k$, respectively, for which (\ref{equalityholds}) holds, then
\begin{equation}\label{RQ}
\bR_n+\bQ_n=\hR_n+\hQ_n
\end{equation}
and the minimal sum $Q_0+\hat{Q}_0$ of the respective variances can serve as a metric quantifying the distance between $\bR$ and $\hR$.

The computation of $\bQ_n,\hQ_n$ minimizing the sum $Q_0+\hat{Q}_0$,  or equivalently minimizing
\[
\frac{1}{n}\trace(\bQ_n+\hQ_n),
\]
is a convex problem --since the positivity constraints are convex. The Toeplitz structure is peripheral, and the idea of defining such metrics extends equally well to non-negative definite Hermitian matrices and to more general positive operators.  For notational convenience we develop the framework in the context of real symmetric matrices.

So, we let
\[
\cM:=\{M\in\mR^{n\times n}\mid M=M^\prime\geq 0\}
\]
be the cone of non-negative symmetric $n\times n$-matrices and
\[
\cT:=\{\bR\in\cM \mid \bR \mbox{ is a Toeplitz matrix}\}
\]
be the cone of non-negative Toeplitz matrices in $\cM$. We address the case of matrices in  $\cM$ and define a suitable metric, which is then specialized to $\cT$.

\begin{proposition}\label{proposition2}{\sf Let $M_1,M_2\in\cM$ and
\begin{eqnarray*}
\tau(M_1,M_2) &:=&\min\left\{\frac{1}{n}\trace(M)\mid M\in\cM,\,\right.\\
&&\left. \phantom{\trace(M)}  M\geq M_1\mbox{ and }M\geq M_2\right\}.
\end{eqnarray*}
Then
\begin{equation}\label{matrixmetric}
\delta(M_1,M_2):=2\tau(M_1,M_2)-\trace(M_1)-\trace(M_2)
\end{equation}
defines a metric on $\cM$.
}\end{proposition}

\begin{proof} Given $M_1,M_2\in\cM$,
\[ \cC(M_1,M_2):=\left\{M\mid  M\geq M_1\mbox{ and }M\geq M_2\right\}
\]
is a (convex) cone of non-negative definite matrices. It follows that there is an element
$M_{12}\in\cC(M_1,M_2)$ having minimal trace.

Clearly $\delta(M_1,M_2)$ is symmetric in its arguments and takes positive values unless $M_1=M_2$, in which case $\delta(M_1,M_2)=0$. Thus, we only need to prove the triangle inequality. Given $M_i\in\cM$ for $i\in\{1,2,3\}$, we denote by 
$M_{ik}$ corresponding minimal elements as above for $i,k\in\{1,2,3\}$, and we let
\[
\Delta_{ik}:=M_{ik}-M_k.
\]
These matrices are non-negative by construction, the identities
\[M_i+\Delta_{ki}=M_k+\Delta_{ik}
\]
hold, and
\[\delta(M_i,M_k)=\frac{1}{n}\trace(\Delta_{ik}+\Delta_{ki})
\]
for $i,k\in\{1,2,3\}$.
But then,
\begin{eqnarray*}
M_1+\Delta_{21}-\Delta_{12}&=&M_2\\
&=&M_3+\Delta_{23}-\Delta_{32},
\end{eqnarray*}
and hence,
\begin{eqnarray*}
M_1+\Delta_{21}+\Delta_{32}&=&M_3+\Delta_{23}+\Delta_{12}.
\end{eqnarray*}
From the minimal property of $\Delta_{13}$ and of $\Delta_{31}$ with regard to having the least value for the combined trace so that $M_1+\Delta_{31}=M_3+\Delta_{13}$, it follows that
\[
\trace(\Delta_{13}+\Delta_{31})\leq \trace(\Delta_{21}+\Delta_{32}+\Delta_{23}+\Delta_{12}).
\]
Therefore,
\begin{eqnarray*}
\delta(M_1,M_2)+\delta(M_2,M_3)&=&\delta(M_1,M_3),
\end{eqnarray*}
which completes the proof.\hfill $\Box$
\end{proof}

We now observe that the steps of the proof of Proposition \ref{proposition2} permit incorporating linear constraints on the structure of elements of $\cM$, such as the constraint of all matrices being Toeplitz. Hence, whereas $\delta(\cdot,\cdot)$ may be used directly as a distance measure between elements of $\cT$, the corresponding minimal-trace perturbations $\Delta_{ik}$ may not belong to $\cT$ in general. But, since the Toeplitz property is a linear constraint, we may define a completely analogous 
distance measure enforcing such perturbations (if so desired) to be Toeplitz.

\begin{proposition}\label{proposition3}{\sf Let $M_1,M_2\in\cT$ and
\begin{eqnarray*}
\tau_{\cTi}(M_1,M_2) &:=&\min\left\{\frac{1}{n}\trace(M)\mid M\in\cT,\mbox{ and } \right.\\
&&\left. \phantom{\frac{1}{n}\mbox{ and also }M\in\cT} M\geq M_1,\, M\geq M_2\right\}.
\end{eqnarray*}
Then
\begin{equation}\label{matrixmetric}
\delta_{\cTi}(M_1,M_2):=2\tau_{\cTi}(M_1,M_2)-\trace(M_1)-\trace(M_2)
\end{equation}
defines a metric on $\cT$.
}\end{proposition}

\begin{proof} The proof follows the steps of the proof of Proposition \ref{proposition2} verbatim, except for the fact that we now constraint all matrices to belong to $\cT$.\hfill $\Box$
\end{proof}

\begin{proposition}\label{proposition4}{\sf Let $f_\by$, $f_\hy$ be power spectral densities,
i.e., nonnegative and integrable on $[-\pi,\,\pi]$. Let as before $\bR_n$, $\hat\bR_n$ denote the  corresponding Toeplitz covariance matrices, and let $n\in\{1,\,2,\,\ldots\}$. Then
\begin{equation}\nonumber
\lim_{n\to\infty}\delta_{\cTi}(\bR_n,\hat\bR_n)=\|f_\by-f_\hy\|_1.
\end{equation}
}\end{proposition}

\begin{proof} Clearly
\begin{equation}\nonumber
\lim_{n\to\infty}\delta_{\cTi}(\bR_n,\hat\bR_n)\leq \|f_\by-f_\hy\|_1
\end{equation}
since a choice of $\psi_k$ and $\hat\psi_k$ with power spectra as in (\ref{psi}-\ref{psihat})
gives rise to partial covariance matrices $\bQ_n$, $\hat\bQ_n$, for all $n$, for which (\ref{RQ}) holds. The respective $0$th elements $Q_0$ and $\hat Q_0$ remain the same for all $n$ and the left hand side is
\[\|f_\by-f_\hy\|_1=Q_0+\hat Q_0
\]
since the power spectra in (\ref{psi}-\ref{psihat}) have no overlap in their support.

To show the converse inequality, consider the sequence of minimizing $\bQ_n$, $\hat\bQ_n$. These are Toeplitz matrices with bounded entries (since their corresponding $0$th element is bounded by $\|f_\by-f_\hy\|_1$). Each can be extended to an infinite Toeplitz matrix, and thereby, gives rise to power spectral densities $q_n$ and $\hat q_n$ such that the first $n$ Fourier coefficients of $f_\by+q_n$ and $f_\hy+\hat q_n$ coincide. The spectral densities $q_n$ and $\hat q_n$ can be obtained from $\bQ_n$, $\hat\bQ_n$ by any particular positive extension, for instance a ``maximum entropy'' one. We can take those as pairs, and since they are bounded there exists a subsequence weakly convergent to  possibly non-negative measures, $d\mu$ and $d\hat\mu$, such that
\[ f_\by d\theta +d\mu = f_\hy d\theta + d\hat\mu
\]
since their Fourier coefficients must coincide.
If $d\mu$, $d\hat\mu$ do have singular parts then these should be identical and the absolutely continuous parts must balance as well, so there exist power spectral densities
$q$ and $\hat q$ such that
\begin{equation}\label{balance} f_\by +q = f_\hy + \hat q.
\end{equation}
But then,
\begin{eqnarray*}
\lim_{n\to\infty}\delta_{\cTi}(\bR_n,\hat\bR_n)
&\geq & \lim_{n\to\infty} \intpi (q_n(\theta)+\hat q_n(\theta))d\theta\\
&=& \intpi (d\mu+d\hat\mu)\\
&\geq& \intpi (q+\hat q)d\theta\\ &\geq& \|f_\by-f_\hy\|_1,
\end{eqnarray*}
the last inequality from (\ref{balance}).
\hfill $\Box$
\end{proof}

\subsection{An example}
The metric $\delta_{\cTi}(\bR_n,\hat\bR_n)$ of the previous section admits no simple expression in terms of the respective eigenvalues. This should be contrasted with its limiting value $d(f_\by,f_\hy)$ which is the $L_1$ distance between the corresponding power spectral densities. We highlight this with an example.

Let
\[
\bR_3=\left[\begin{matrix}\phantom{1}1\phantom{1} & \phantom{1}1\phantom{1} & \phantom{1}1\phantom{1}\\ 1 & 1 & 1\\ 1& 1& 1\end{matrix}\right]
\]
and
\[
\hat\bR_3=\left[\begin{matrix}1 & 1/2 & 1/2\\ 1/2 & 1 & 1/2\\ 1/2& 1/2& 1\end{matrix}\right].
\]
Then, clearly,
\[
\bQ_3=\left[\begin{matrix}\phantom{1}x\phantom{1} & \phantom{1}y\phantom{1} & \phantom{1}y\phantom{1}\\ y & x & y\\ y& y& x\end{matrix}\right]
\]
and
\[
\hat\bQ_3=\left[\begin{matrix}\phantom{1}x\phantom{1} & \phantom{1}v\phantom{1} & \phantom{1}v\phantom{1}\\ v & x & v\\ v& v& x\end{matrix}\right]
\]
where
\[
1+y=1/2+v
\]
and $x$ is minimal subject to $\bQ_3\geq 0$ as well as $\hat\bQ_3\geq 0$. The last two inequalities imply that
\begin{eqnarray*}
0\leq &x & \mbox{ as well as }\\
-\frac12 x\leq &y,v & \leq x.
\end{eqnarray*}
It follows that the optimal choice (minimal $x$) is
\begin{eqnarray*}
x &=& \phantom{-}1/3\\
y&=& -1/6\\
v&=& \phantom{-}1/3.
\end{eqnarray*}
Then
\[\delta_{\cTi}(\bR_3,\hat\bR_3) = 2x= 2/3,
\]
while the respective eigenvalues are
\begin{eqnarray*}
{\rm spec}(\bR_3)&=&\{3,\,0,\,0\}\mbox{ and}\\
{\rm spec}(\hat\bR_3)&=&\{1,\,1,\,1/2\}.
\end{eqnarray*}
It appears that there is no simple expression for $\delta_{\cTi}(\bR_3,\hat\bR_3)$ based solely on knowledge of ${\rm spec}(\bR_3)$ and ${\rm spec}(\hat\bR_3)$.

The covariance $\bR_3$ has a unique extension and corresponds to a measure with unit weight at $\theta=0$, i.e., a spectral line (Dirac delta) at $\theta=0$. Assuming that $\hat \bR_3$ originates from a spectral measure which has a similar weight of amplitude $1/2$ at $\theta=0$ and a uniform absolutely continuous part of amplitude $1/2$, then
\[
\|d\mu-d\hat\mu\|_1=1/2 + 1/2 =1
\]
adding the $L_1$-norm of the difference of the absolutely continuous parts with the absolute integral of the discrepancy between the two measures. We leave it as an exercise to the reader to verify that if $\hat\bR_n$ is as we just assumed, namely $\hat R_k=1/2$ for $k\geq 1$, and similarly, $R_k=1$ for all $k$, then $\delta_{\cTi}(\bR_n,\hat\bR_n)\to 1$ as $n\to \infty$.

\section{Approximating sample covariances}

It is often the case that the autocovariance matrix $\bR_n$ of a random process $\by_k$ is estimated in a way that does not guarantee this to be Toeplitz. For instance, it is quite common for $\bR_n$ to be estimated by averaging observation samples
\[
\hat\bR_{n}= \frac{1}{N+1}\sum_{\ell=0}^N \left[\begin{array}{c}y_{1+\ell}\\\vdots\\y_{n+\ell}\end{array}\right]
\left[\begin{array}{ccc}y_{1+\ell}&\ldots&y_{n+\ell}\end{array}\right]
\]
The estimate $\hat\bR_{n}$ is non-negative definite by construction, but may not be Toeplitz. Yet, for purposes of analysis it is often beneficial to approximate $\bR_{n}$ by a Toeplitz one, or possibly, by one with additional structure (e.g., corresponding to a moving average process or, more generally, to the state of a known dynamical system). The problem of seeking such an approximant which is closest to $\bR_{n}$ in $\delta(\cdot,\cdot)$, is readily solvable via convex optimization.

\subsection{Comparison with the von Neumann entropy}
In \cite{kimurafest}, the question was raised as to what are appropriate ways to approximate a given sample covariance with one that abides by a known linear structure.
It was proposed that the Kullback-Leibler-von Neuman distance
\[{\mathbb{S}}(\hat \bR \|\bR):={\rm trace}\left(\hat \bR\left(\log \hat \bR -\log \bR\right)\right)
\]
provides a convenient convex functional for which the optimal approximant is uniquely defined. An academic example was presented in \cite{kimurafest} which is recapitulated here as it helps underscore differences with approximation in the sense of minimizing $\delta(\hat\bR,\bR)$.

Consider the positive-definite matrix below as the estimated value for a covariance matrix
\[
\hat\bR_3=\frac{1}{3}\left[\begin{matrix}1.1 & .9 & 1.05\\
                                                          .9 & .8 & .9\\
                                                          1.05 & .9 & 1.1\end{matrix}\right].
\]
The minimizer of
\[
\{{\mathbb{S}}(\hat\bR,\bR)\mid \bR \mbox{ being Toeplitz},\, \bR>0,\,{\rm trace}(\bR)={\rm trace}(\hat\bR)\}
\]
is unique (see \cite{kimurafest}) and given by
\[
\bR_{3,\rm vN}=\frac{1}{3}\left[\begin{matrix}1 & .942 & .957\\
                                                          .942 & 1 & .942\\
                                                          .957 & .942 & 1\end{matrix}\right].
\]                                                       
It is interesting to point out the the closest Toeplitz matrix to $\hat\bR$ in the least-squares sense fails to be positive-definite (\cite{kimurafest}, cf.\ \cite{IT}).
On the other hand, the optimal approximant in $\delta(\cdot,\cdot)$-sense can be obtained by observation and is equal to
\[
\bR_{3,\delta}=\frac{1}{3}\left[\begin{matrix}1.1 & .9 & 1.05\\
                                                          .9 & 1.1 & .9\\
                                                          1.05 & .9 & 1.1\end{matrix}\right].
\]
In the above, a second subscript indicates the sense in which the matrix approximates $\hat\bR_3$.
Obviously the traces of $\bR_{3,\delta}$ and $\hat\bR_3$ are not the same, in general. However, equality of the traces can be easily imposed as an added linear constraint.

\subsection{Structured covariances}

For purposes of illustration, consider a moving average process
\[ \by_k=\bw_k+\bw_{k-1}+\bw_{k-2}
\]
where $\bw_k$ is a zero-mean, unit-variance, Gaussian white noise process. The autocorrelation sequence of $\by_k$ is
\[
\left[\begin{matrix} R_0&R_1 &R_2 &R_3& 0 &\ldots\end{matrix}\right]=
\left[\begin{matrix} 3&2 &1 &0& 0 &\ldots\end{matrix}\right].
\]
Simulating $\by_k$ over a window $k\in\{0,1,\ldots,100\}$, and based on a particular such realization, the corresponding $n\times n$ sample covariance matrix, for $n=5$, was computed to be
\[
\hat\bR_5=\left[\begin{matrix}
    4.0362   & 2.9053    &1.8043   & 0.4042    &0.1718\\
    2.9053    &4.0547    &2.9268    &1.7945    &0.3800\\
    1.8043    &2.9268    &4.0792    &2.9143    &1.7733\\
    0.4042    &1.7945    &2.9143    &4.0819    &2.9421\\
    0.1718    &0.3800    &1.7733    &2.9421   & 4.0237
\end{matrix}\right].
\]
Obviously, this matrix is not Toeplitz due to the finiteness of the observation record.
The closest Toeplitz approximant to $\hat\bR_5$, in the sense of the metric $\delta(\cdot,\cdot)$, turns out to be
\[
\bR_{5,\rm Toeplitz}=\left[\begin{matrix}
4.0677   & 2.9237    &1.7912    &0.3979   & 0.1822\\
    2.9237  &  4.0677   & 2.9237   & 1.7912  &  0.3979\\
    1.7912  &  2.9237   & 4.0677   & 2.9237  &  1.7912\\
    0.3979  &  1.7912    &2.9237   & 4.0677   & 2.9237\\
    0.1822  &  0.3979    &1.7912   & 2.9237  &  4.0677\end{matrix}\right]
\]
for which
\[
\delta(\hat\bR_5,\bR_{5,\rm Toeplitz})=0.0308.
\]
Interestingly, $\bR_{5,\rm Toeplitz}$ does not correspond to a moving average process of order $2$ (or even, of order $3,4$) as it can be readily verified by the fact that the trigonometric polynomials, e.g.,
\[\sum_{k=-4}^4 R_k e^{jk\theta}
\]
takes negative values.

The set of covariance matrices which are generated by moving average processes of a given order, is convex and admits a characterization via a set of linear matrix inequalities (\cite{StoicaMA,GeoDecomposition}). Thus, the closest approximant to $\hat\bR$ which corresponds to a moving average process of any given order can be readily computed. In particular, if we specify the order to be $2$, then the optimal approximant to $\hat \bR_5$ becomes
\[
\bR_{5,\rm MA(2)}=\left[\begin{matrix}
    3.9945&    2.1588&    0.5693  &       0      &   0\\
    2.1588 &   3.9945  &  2.1588   & 0.5693   &      0\\
    0.5693 &   2.1588    &3.9945   & 2.1588 &   0.5693\\
         0    &0.5693    &2.1588   & 3.9945 &   2.1588\\
         0     &    0    &0.5693   & 2.1588 &   3.9945
\end{matrix}\right]
\]
for which
\[
\delta(\hat\bR_5,\bR_{5,\rm MA(2)})=
  1.2161.
\]

\printindex
\end{document}